\documentclass{article}
\usepackage[english]{babel}
\usepackage[cp1251]{inputenc}
\usepackage{amsmath,amsthm,amssymb,amstext,amscd,amsfonts}
\begin{document}
\begin{center}
\textbf{\Large{Double-Layer Potentials for a Generalized Bi-Axially
Symmetric Helmholtz Equation}}\\

\medskip
\textbf{Junesang Choi${}^1$, Anvar Hasanov${}^2$ and Hari M. Srivastava$^{3,\ast}$}\\
\medskip
${}^1$ Department of Mathematics, Dongguk University,\\
 Gyeongju 780-714, Republic of Korea\\
\textbf{E-Mail: junesang@mail.dongguk.ac.kr}\\
${}^2$Department of Mathematics, I. M. Gubkin Russian State University of Oil and Gas, \\
Tashkent Branch,  Tashkent 100180, Uzbekistan\\
\textbf{E-mail: anvarhasanov@yahoo.com}\\
${}^3$ Department of Mathematics and Statistics, University of
Victoria,\\ Victoria, British Columbia V8W 3P4, Canada.\\
 \textbf{E-Mail: harimsri@math.uvic.ca}\\

\textbf{$^{\ast}$Corresponding Author}\\

\end{center}
\medskip

\begin{abstract}

The double-layer potential plays an important r$\hat{\rm o}$le in solving boundary value problems of elliptic equations.
Here, in this paper,  we aim at  introducing and investigating double layer potentials for a
generalized bi-axially symmetric Helmholtz equation. By using some properties of one of  Appell's hypergeometric functions
in two variables,  we prove limiting theorems  and derive integral equations concerning a denseness of double-layer potentials.
\end{abstract}



\vskip 15 mm
\section*{\bf 1. Introduction}

Potential theory has played a paramount r$\hat{\rm o}$le in both analysis and computation for boundary value problems for elliptic partial differential equations. Numerous applications can be found in fracture mechanics, fluid mechanics, elastodynamics, electromagnetics, and acoustics. Results from potential theory allow us to represent boundary value problems in integral equation form. For problems with known Green's functions, an integral equation formulation leads to powerful numerical approximation schemes.

The double-layer potential plays an important r$\hat{\rm o}$le in solving boundary value problems of elliptic equations. The representation of the solution of the (first) boundary value problem is sought as a double-layer potential with unknown density  and an application of certain property leads to a Fredholm equation of the second kind for determining the function (see \cite{gunt-67} and \cite{miran-70}).

By applying a method of complex
analysis (based upon analytic functions), Gilbert \cite{Gilb-69}
constructed an integral representation of solutions of the following  generalized
bi-axially Helmholtz equation:
$$
H_{\alpha ,\beta }^\lambda \left( u \right) \equiv u_{xx}  +
u_{yy}  + \displaystyle \frac{{2\alpha }}{x}u_x  + \displaystyle
\frac{{2\beta }}{y}u_y  - \lambda ^2 u = 0  \quad \left(0 <\alpha,\, \beta  < \frac{1}{2}\right),\eqno \left( {H_{\alpha
,\beta }^\lambda  } \right)
$$
where $\alpha$, $\beta$ and $\lambda$ are
constants. When $\lambda  = 0$, this equation is known as the equation of
the generalized axially symmetric potential theory  whose name
is due to Weinstein who first considered fractional dimensional
space in potential theory (see \cite{Wein-48} and \cite{Wein-53}). The special case where $
\,\lambda  = 0$ was also investigated by (among others) Erd\'{e}lyi (see \cite{Erde-56} \cite{Erde-65}),
Gilbert (see \cite{Gilb-60},   \cite{Gilb-64},  \cite{Gilb-65},  \cite{Gilb-67} and \cite{Gilb-68}),  Gilbert and Howard \cite{Gi-Ho}  , Ranger  \cite{Rang} and  Henrici (see \cite{Henr-53} and \cite{Henr-60}).
 Certain problems for the equation $H_{\alpha ,\beta }^\lambda$ were studied by many authors (see, for example,
 \cite{Alti-1},   \cite{Alti-2},  \cite{Al-Eu},  \cite{Frya},  \cite{Hube},  \cite{Kuma}, \cite{Lo},  \cite{Mari}, \cite{McCo-79},
   \cite{McCo-80},  \cite{Pi-Bo},   \cite{Ra-Ha} and \cite{Wein}).

 Fundamental solutions of the equation $\left(H_{\alpha ,\beta }^\lambda\right)$ are
constructed recently (see \cite{Hasa-07}). In fact, the fundamental solutions of the equation $\left(H_{\alpha ,\beta }^\lambda\right)$
when $\,\lambda  = 0$ can be expressed in terms of   Appell's hypergeometric function in two variables of the second kind, that is,
the Appell function
$$F_2 \left( {a,b_1 ,b_2 ;c_1 ,c_2 ;x,y} \right)$$
 defined by (see \cite[p. 224, Eq. 5.7.1 (7)]{Er-Ma-Ob-Tr}; see also \cite[p. 14, Eq.~(12)]{Ap-Ka} and \cite[p. 23, Eq. 1.3 (3)]{Sr-Ka})
 $$
F_2 \left( {a,b_1 ,b_2 ;c_1 ,c_2 ;x,y} \right) = \displaystyle
\sum\limits_{m,n = 0}^\infty  {} \displaystyle \frac{{\left( a
\right)_{m + n} \left( {b_1 } \right)_m \left( {b_2 } \right)_n
}}{{\left( {c_1 } \right)_m \left( {c_2 } \right)_n}}\,\frac{x^m}{m!} \,\frac{y^n}{n!},
\eqno (1.1)
$$
where $(\lambda)_\nu$ denotes the Pochhammer symbol or the \emph{shifted factorial}, since
$$ (1)_n =n! \quad \left(n \in \mathbb{N}_0:=\mathbb{N}\cup \{0\}; \,\,\mathbb{N}:=\{1,\,2,\,3,\,\cdots\}\right),$$
which is defined (for $\lambda,\,\nu \in \mathbb{C}$), in terms of the familiar Gamma function, by
$$  (\lambda)_\nu := \frac{\Gamma (\lambda +\nu)}{ \Gamma (\lambda)}
 =\left\{\aligned & 1  \hskip 41 mm (\nu=0;\,\, \lambda \in \mathbb{C}\setminus \{0\}) \\
        & \lambda (\lambda +1) \cdots (\lambda+n-1) \hskip 8mm (\nu=n \in {\mathbb N};\,\, \lambda \in \mathbb{C}),
   \endaligned \right. \\
$$
it being \emph{understood conventionally} that $(0)_0:=1$. We thus obtain the following results:
$$
q_1 \left( {x,y;x_0 ,y_0 } \right) = \displaystyle k_1 \left( {r^2
} \right)^{ - \alpha  - \beta } F_2 \left( {\alpha  + \beta, \alpha ,\beta ;2\alpha ,2\beta ;\xi ,\eta } \right),\eqno (1.2)
$$
$$
q_2 \left( {x,y;x_0 ,y_0 } \right) = \displaystyle k_2 \left( {r^2
} \right)^{\alpha  - \beta  - 1} x^{1 - 2\alpha } x_0^{1 - 2\alpha
} F_2 \left( {1 - \alpha  + \beta,1 - \alpha ,\beta ;2 - 2\alpha
,2\beta ;\xi ,\eta } \right),\eqno (1.3)
$$
$$
q_3 \left( {x,y;x_0 ,y_0 } \right) = \displaystyle k_3 \left( {r^2
} \right)^{ - \alpha  + \beta  - 1} y^{1 - 2\beta } y_0^{1 -
2\beta } F_2 \left( {1 + \alpha  - \beta ,\alpha ,1 - \beta
;2\alpha ,2 - 2\beta ;\xi ,\eta } \right)\eqno (1.4)
$$
and
\begin{equation}\tag{1.5} \aligned &
q_4 \left( {x,y;x_0 ,y_0 } \right) = \displaystyle k_4 \left( {r^2
} \right)^{\alpha  + \beta  - 2} x^{1 - 2\alpha } y^{1 - 2\beta }
x_0^{1 - 2\alpha } y_0^{1 - 2\beta } \\
& \hskip 25mm \cdot  F_2 \left( {2 - \alpha  -
\beta,1 - \alpha ,1 - \beta ;2 - 2\alpha ,2 - 2\beta ;\xi ,\eta }
\right),\endaligned
\end{equation}
where
$$
k_1  =\displaystyle \frac{{2^{2\alpha  + 2\beta } }}{{4\pi
}}\frac{{\Gamma \left( \alpha  \right)\Gamma \left( \beta
\right)\Gamma \left( {\alpha  + \beta } \right)}}{{\Gamma \left(
{2\alpha } \right)\Gamma \left( {2\beta } \right)}},\eqno (1.6)
$$
$$
k_2  = \displaystyle \frac{{2^{2 - 2\alpha  + 2\beta } }}{{4\pi
}}\frac{{\Gamma \left( {1 - \alpha } \right)\Gamma \left( \beta
\right)\Gamma \left( {1 - \alpha  + \beta } \right)}}{{\Gamma
\left( {2 - 2\alpha } \right)\Gamma \left( {2\beta }
\right)}},\eqno (1.7)
$$
$$
k_3  = \displaystyle \frac{{2^{2 + 2\alpha  - 2\beta } }}{{4\pi
}}\frac{{\Gamma \left( \alpha  \right)\Gamma \left( {1 - \beta }
\right)\Gamma \left( {1 + \alpha  - \beta } \right)}}{{\Gamma
\left( {2\alpha } \right)\Gamma \left( {2 - 2\beta }
\right)}},\eqno (1.8)
$$
$$
k_4  = \displaystyle \frac{{2^{4 - 2\alpha  - 2\beta } }}{{4\pi
}}\frac{{\Gamma \left( {1 - \alpha } \right)\Gamma \left( {1 -
\beta } \right)\Gamma \left( {2 - \alpha  - \beta }
\right)}}{{\Gamma \left( {2 - 2\alpha } \right)\Gamma \left( {2 -
2\beta } \right)}},\eqno (1.9)
$$

$$
\left. {\begin{array}{*{20}c}
   {r^2 }  \\
   {r_1^2 }  \\
   {r_2^2 }  \\
\end{array}} \right\} = \left( {x\begin{array}{*{20}c}
    -   \\
    +   \\
    -   \\
\end{array}x_0 } \right)^2  + \left( {y\begin{array}{*{20}c}
    -   \\
    -   \\
    +   \\
\end{array}y_0 \,} \right)^2 ,
\,\,\, \xi  = \displaystyle \frac{{r^2  - r_1^2 }}{{r^2 }} \,\,\,\text{and}\,\,\,
\eta  = \displaystyle \frac{{r^2 - r_2^2 }}{{r^2 }}.\eqno (1.10)
$$

The fundamental solutions given by (1.2) to (1.5) possess
the following properties:
$$
\left. {x^{2\alpha }} \displaystyle \frac{\partial}{\partial x}
\, \left\{q_1 \left({x,y;x_0 ,y_0 } \right)  \right\}\right|_{x = 0}  = 0,
\quad \left. {y^{2\beta }}\frac{\partial}{\partial y} \left\{ q_1 \left( {x,y;x_0 ,y_0
} \right)\right\} \right|_{y = 0}  = 0,\eqno (1.11)
$$

$$
\left. q_2 \left( x,y;x_0 ,y_0  \right) \right|_{x = 0}  = 0,
\quad \left. {y^{2\beta }} \displaystyle
\frac{\partial}{\partial y} \left\{q_2
\left( x,y;x_0 ,y_0 \right)\right\}\right|_{y = 0}  =
0,\eqno (1.12)
$$

$$
\left. {x^{2\alpha }} \displaystyle \frac{\partial}{\partial x}  \left\{ q_3 \left(
{x,y;x_0 ,y_0 } \right) \right\} \right|_{x = 0}  =
0,\quad \left. q_3 \left( {x,y;x_0 ,y_0 } \right) \right|_{y =
0}  = 0,\eqno (1.13)
$$
$$
\left. {q_4 \left( {x,y;x_0 ,y_0 } \right)} \right|_{x = 0}  = 0
\quad \text{and} \quad \left. {q_4 \left( {x,y;x_0 ,y_0 } \right)} \right|_{y = 0}
= 0.\eqno (1.14)
$$

Here, by making use of the fundamental solutions given by (1.2) to (1.5) in the following domain:
\begin{equation}\tag{1.15}
 \Omega  \subset \mathbb{R}_ + ^2  = \left\{ {\left( {x,y}
\right)\, :\,x > 0 \quad \text{and} \quad  y > 0} \right\},
\end{equation}
  we aim at investigating a double-layer potential for
the equation $\left(H_{\alpha ,\beta }^0\right)$. Furthermore, we prove some results (see Lemmas 1 to 3 and Theorem 3) on limiting values of the double-layer potential in (3.2). These results are (potentially) useful for future works in which boundary value problems for the equation
 $\left(H_{\alpha ,\beta }^0\right)$
are investigated in more general domains.

\vskip 5 mm
\section*{\bf 2. Green's formula}

We begin by considering the following identity:
$$
\displaystyle x^{2\alpha } y^{2\beta } \left[ {uH_{\alpha ,\beta
}^0 \left( v \right) - vH_{\alpha ,\beta }^0 \left( u \right)}
\right] = \displaystyle \frac{\partial }{{\partial x}}\left[
{x^{2\alpha } y^{2\beta } \left( {v_x u - vu_x } \right)} \right]
+ \displaystyle \frac{\partial }{{\partial y}}\left[ {x^{2\alpha }
y^{2\beta } \left( {v_y u - vu_y } \right)} \right].\eqno (2.1)
$$
Integrating both parts of the identity (2.1) on a domain $\Omega$ in (1.15), and using Green's formula, we find that
$$
\displaystyle  \int\int_\Omega\, x^{2\alpha } y^{2\beta } \left[ {uH_{\alpha ,\beta
}^0 \left( v \right) - vH_{\alpha ,\beta }^0 \left( u \right)}
\right]\,dxdy = \displaystyle \int\limits_S \, x^{2\alpha }
y^{2\beta } u\left( {v_x dy - v_y dx} \right) - x^{2\alpha }
y^{2\beta } v\left( {u_x dy - u_y dx} \right),\eqno (2.2)
$$
where $S = \partial \Omega $ is the boundary of the domain
$\Omega$.

 If $ u\left( {x,y} \right)$ and $ v\left( {x,y} \right)$
are solutions of the equation $\left(H_{\alpha ,\beta }^0 \right)$, we find from
 (2.2) that
$$
\displaystyle \int\limits_S {} x^{2\alpha } y^{2\beta } \left( {u
\displaystyle \frac{{\partial v}}{{\partial n}} - v \displaystyle
\frac{{\partial u}}{{\partial n}}} \right)ds = 0,\eqno (2.3)
$$
where
$$
\displaystyle \frac{\partial }{{\partial n}} = \displaystyle
\frac{{dy}}{{ds}} \displaystyle \frac{\partial }{{\partial x}} -
\frac{{dx}}{{ds}} \displaystyle \frac{\partial }{{\partial
y}},\,\,\,\displaystyle \frac{{dy}}{{ds}} = \cos \left( {n,x}
\right) \quad \text{and} \quad \displaystyle \frac{{dx}}{{ds}} =  - \cos \left(
{n,y} \right),\eqno (2.4)
$$
$n$ being the exterior normal to the curve $S$. We also obtain the following identity:
$$ \int\int_\Omega\,
\displaystyle x^{2\alpha } y^{2\beta } \left[ {u_x^2  + u_y^2 }
\right]\, dxdy =\displaystyle \int\limits_S {} x^{2\alpha } y^{2\beta
} u \displaystyle \frac{{\partial u}}{{\partial n}}ds,\eqno (2.5)
$$
where $ u\left( {x,y} \right)$ is a solution of the equation
$\left(H_{\alpha ,\beta }^0\right)$. The special case of (2.3) when  $v =
1$ reduces to the following form:
$$
\displaystyle \int\limits_S {} x^{2\alpha } y^{2\beta }
\displaystyle \frac{{\partial u}}{{\partial n}}ds = 0.\eqno (2.6)
$$
We note from (2.6) that the integral of the normal derivative of a solution of the equation
$\left(H_{\alpha ,\beta }^0\right) $ with a weight $x^{2\alpha } y^{2\beta }$
along the boundary $S$ of the domain $\Omega $ in (1.15) is equal to zero.

\vskip 15 mm
\section*{\bf 3. A Double-Layer Potential $w^{\left( 1 \right)} \left( {x_0 ,y_0 } \right)$}

Let $\Omega$ in (1.15) be a domain bounded by intervals $\left( {0,a}
\right)$ and $\left( {0,b} \right)$ of the $x-$ and $y-$ axes,
respectively, and a curve $\Gamma $ with the extremities at points
$A\left( {a,0} \right)$ and $B\left( {0,b} \right)$. The parametric
equation of the curve $\Gamma $ is given by
 $$x = x\left( s \right) \quad \text{and} \text    y =
y\left( s \right),$$
 where $s$ denotes  the length of an arc beginning from the
point $A\left( {a,0} \right)$. We
 assume the following properties of the curve $\Gamma $:
 \begin{enumerate}
   \item[(i)] The functions $x = x\left( s \right)$ and $y =y\left( s \right)$ have
continuous derivatives $x'\left( s \right)$ and $y'\left( s \right)$ on
a segment $\left[ {0,l} \right]$ and do not vanish simultaneously;

   \item[(ii)] The second derivatives $x''\left( s \right)$ and $y''\left( s \right)$ satisfy
the  H$\ddot{\rm o}$lder condition on $\left[ {0,l} \right]$,
where $l$  denotes the length of the curve $\Gamma $;

   \item[(iii)] In some neighborhoods of points $A\left( {a,0} \right)$ and
$B\left( {0,b} \right)$, the following conditions are satisfied:
$$
\left| {\displaystyle \frac{{dx}}{{ds}}} \right| \le cy^{1 +
\varepsilon } \left( s \right) \quad \text{and} \quad \left| {\displaystyle
\frac{{dy}}{{ds}}} \right| \le cx^{1 + \varepsilon } \left( s
\right)  \quad (0 < \varepsilon  < 1;\,\,c = \text{a constant}), \eqno (3.1)
$$
  $\left( {x,y} \right)$ being the coordinates
of a variable point on the curve $\Gamma $.
 \end{enumerate}

  Consider the following integral
$$
w^{\left( 1 \right)} \left( {x_0 ,y_0 } \right) = \displaystyle
\int_0^l \,{x^{2\alpha } y^{2\beta } } \mu _1 \left( s
\right) \displaystyle \frac{\partial}{\partial n}\,\left\{ q_1 \left( {x,y;x_0 ,y_0 }
\right)\right\}\,ds,\eqno (3.2)
$$
where $\mu _1 \left( s \right) \in C\left[ {0,\,l} \right]$
and $q_1$ is given in (1.2). We call the integral (3.2) a \emph{double-layer potential with
denseness} $\mu _1 \left( s \right)$.

We now investigate some properties of a double-layer potential  $w^{\left( 1 \right)} \left( {x_0 ,y_0 } \right)$ with
denseness $\mu _1 \left( s \right)$.

\vskip 5mm
\textbf{Lemma 1.} \emph{The following formula holds true}:
$$
w_1^{\left( 1 \right)} \left( {x_0 ,y_0 } \right) =
\left\{
  \begin{array}{rl}
  - 1 & \qquad  \left(\left( {x_0 ,y_0 } \right) \in \Omega\right) \\
    - \frac{1}{2} & \qquad   \left(\left( {x_0 ,y_0 } \right) \in \Gamma\right) \\
    0 & \qquad  \left(\left( {x_0 ,y_0 } \right) \not\in \bar \Omega\right),
  \end{array}
\right.\eqno (3.3)
$$
\emph{where a domain} $\Omega$ \emph{and the curve} $\Gamma$ \emph{are described as in this section and}
$\bar \Omega:=\Omega \cup \Gamma$.

\begin{proof}
    \textbf{Case 1.} When $\left( {x_0 ,y_0 } \right) \in \Omega$,
 we cut a circle    centered at $\left( {x_0 ,y_0 } \right)$  with a small radius
$\rho$  off  the domain $\Omega$
and denote the remaining part by $\Omega ^\rho$ and  the circuit of the cut-off-circle by
$C_\rho$.
 The function $q_1 \left( {x,y;x_0 ,y_0 } \right)$  in (1.2) is a regular solution
of the equation $\left(H_{\alpha ,\beta }^0\right) $ in the domain $\Omega
^\rho$. Using the following derivative formula of Appell's hypergeometric function $F_2$ (see \cite[p. 19, Eq.~(20)]{Ap-Ka}):
$$\aligned
& \frac{\partial ^{m + n}}{\partial x^m \partial y^n}\,\left\{ F_2 \left( a,b_1 ,b_2 ;c_1
,c_2 ;x,y \right)\right\}\\
& \hskip 5mm  =
\frac{{\left( a \right)_{m + n} \left( {b_1 } \right)_m \left(
{b_2 } \right)_n }}{{\left( {c_1 } \right)_m \left( {c_2 }
\right)_n }}F_2 \left( {a + m + n,b_1  + m,b_2  + n;c_1  + m,c_2
+ n;x,y} \right),\endaligned  \eqno (3.4)
$$
we have
$$
\begin{array}{l}
\displaystyle  \frac{\partial}{\partial x}\, \left\{ q_1 \left( {x,y;x_0 ,y_0 }
\right)\right\} =  - 2\left( {\alpha  + \beta } \right)k_1
\left( {r^2 } \right)^{ - \alpha  - \beta  - 1} \left( {x - x_0 }
\right)F_2 \left( {\alpha  + \beta ,\alpha ,\beta ;2\alpha ,2\beta ;\xi ,\eta } \right)  \\
 \hskip 3mm - 2\left( {\alpha  + \beta } \right)k_1 \left( {r^2 } \right)^{ - \alpha  -
\beta  - 1} x_0 F_2 \left( {\alpha  + \beta  + 1,\alpha  + 1,
\beta ;2\alpha  + 1,2\beta ;\xi ,\eta } \right)   \\
 \hskip 3mm - 2k_1 \left( {r^2 } \right)^{ - \alpha  - \beta  - 1} \left( {x -
x_0 } \right)\left[ {\displaystyle \frac{{\left( {\alpha  + \beta
} \right)\alpha }}{{2\alpha }}\xi F_2 \left( {\alpha  +
\beta  + 1,\alpha  + 1,\beta ;2\alpha  + 1,2\beta ;\xi ,\eta } \right)} \right.   \\
\left.{\,\,\,\,\,\,\,\,\,\,\,\,\,\,\,\,\,\,\,\,\,\,\,\,\,\,\,\,\,\,\,\,\,
\,\,\,\,\,\,\,\,\,\,\,\,\,\,\,\,\,\,\,\,\,\,\,\,\,\,\,\,\,\,\,\,\,\,  +
\displaystyle \frac{{\left( {\alpha + \beta } \right)\beta
}}{{2\beta }}\eta F_2\left( {\alpha  + \beta  + 1,\alpha ,
\beta  + 1;2\alpha ,2\beta  + 1;\xi ,\eta } \right)} \right]. \\
\end{array}\eqno
(3.5)
$$
By applying the following known contiguous relation (see \cite[p. 21]{Ap-Ka}):
$$
\begin{array}{l}
\displaystyle \frac{{b_1 }}{{c_1 }}xF_2 \left( {a + 1,b_1  + 1,b_2
;c_1  + 1,c_2 ;x,y}\right) + \displaystyle \frac{{b_2 }}{{c_2 }}yF_2
\left( {a + 1,b_1 ,b_2  + 1;c_1 ,c_2  + 1;x,y} \right) \\
 \hskip 10mm = F_2 \left( {a + 1,b_1 ,b_2 ;c_1 ,c_2 ;x,y} \right) - F_2 \left( {a,b_1 ,b_2 ;c_1 ,c_2 ;x,y} \right) \\
\end{array}\eqno (3.6)
$$
to (3.5), we obtain
$$\aligned
 \frac{\partial}{\partial x} \left\{ q_1 \left( {x,y;x_0 ,y_0 }
\right)\right\}  =&  - 2\left( {\alpha  + \beta } \right)k_1
x_0 \left( {r^2 } \right)^{ - \alpha  - \beta  - 1} F_2 \left(
{\alpha  + \beta  + 1,\alpha  + 1,\beta ;2\alpha  + 1,2\beta ;\xi ,\eta } \right) \\
&  - 2\left( {\alpha  + \beta } \right)k_1 \left( {x - x_0 } \right)\left( {r^2 } \right)^{ -
\alpha  - \beta  - 1} F_2 \left( {\alpha  + \beta  + 1,\alpha ,\beta ;2\alpha ,2\beta ;\xi ,\eta } \right).
\endaligned \eqno (3.7)
$$
Similarly, we find that
$$
\aligned   \frac{\partial}{\partial y} \left\{ q_1 \left( {x,y;x_0 ,y_0 }
\right)\right\} = &  - 2\left( {\alpha  + \beta } \right)k_1
y_0 \left( {r^2 } \right)^{ - \alpha  - \beta  - 1} F_2 \left(
{\alpha  + \beta  + 1,\alpha ,1 + \beta ;2\alpha ,1 + 2\beta ;\xi ,\eta } \right) \\
& - 2\left( {\alpha  + \beta } \right)k_1 \left( {y - y_0 } \right)\left( {r^2 } \right)^{ -
\alpha  - \beta  - 1} F_2 \left( {\alpha  + \beta  + 1,\alpha ,\beta ;2\alpha ,2\beta ;\xi ,\eta } \right). \\
\endaligned \eqno (3.8)
$$
Thus, with the help of (3.7) and (3.8),  it follows from (1.2) and (2.4) that
$$
\begin{array}{l}
\displaystyle \frac{\partial}{\partial n} \left\{ q_1 \left( {x,y;x_0 ,y_0 }
\right)\right\} =  - \left( {\alpha  + \beta } \right)k_1
\left( {r^2 } \right)^{ - \alpha  - \beta } F_2 \left( {\alpha  +
\beta  + 1,\alpha ,\beta ;2\alpha ,2\beta ;\xi ,\eta }
\right)\displaystyle \frac{\partial }{{\partial n}}\left\{ \ln \,r^2  \right\}\\
\hskip 15mm + 2\left( {\alpha  + \beta } \right)k_1 y_0 \left( {r^2 }
\right)^{ - \alpha  - \beta  - 1} F_2 \left( {\alpha  + \beta  +
1,\alpha ,1 + \beta ;2\alpha ,1 + 2\beta ;\xi ,\eta }
\right) \displaystyle \frac{d}{ds} \left\{ x\left( s \right)\right\} \\
\hskip 15mm - 2\left( {\alpha  + \beta } \right)k_1 x_0 \left( {r^2 }
\right)^{ - \alpha  - \beta  - 1} F_2 \left( {\alpha  + \beta  +
1,\alpha  + 1,\beta ;2\alpha  + 1,2\beta ;\xi ,\eta }
\right)\displaystyle \frac{d}{ds} \left\{ y\left( s \right)\right\}. \\
\end{array}\eqno(3.9)
$$
Applying (2.6) and considering the identity (1.11),  we get
the following formula:
$$
w_1^{\left( 1 \right)} \left( {x_0 ,y_0 } \right) = \displaystyle
\mathop {\lim }\limits_{\rho  \to 0} \int\limits_{C_\rho  } {}
x^{2\alpha } y^{2\beta } \displaystyle \frac{\partial}{\partial n} \left\{ q_1 \left(
{x,y;x_0 ,y_0 } \right)\right\}ds.\eqno (3.10)
$$
Substituting from (3.9) into (3.10), we find that
$$
\aligned
& w_1^{\left( 1 \right)} \left( {x_0 ,y_0 } \right)\\
& \hskip 3mm =  - \left( {\alpha  + \beta } \right)k_1
\displaystyle \mathop
{\displaystyle \lim }\limits_{\rho \to 0} \displaystyle
\int\limits_{C_\rho  } {} x^{2\alpha } y^{2\beta } \left( {r^2 } \right)^{ - \alpha  - \beta } F_2 \left( {\alpha  +
\beta  + 1,\alpha ,\beta ;2\alpha ,2\beta ;\xi ,\eta }
\right)\displaystyle \frac{\partial }{{\partial n}}\left\{ \ln \,r^2  \right\} ds\\
& \hskip 6mm - 2\left( \alpha  + \beta  \right)k_1 x_0 \mathop {\lim }
\limits_{\rho \to 0} \displaystyle \int\limits_{C_\rho  } {}
x^{2\alpha } y^{2\beta } \left( {r^2 } \right)^{ - \alpha  - \beta
- 1} F_2 \left( {1 + \alpha  + \beta ,1 + \alpha ,\beta ;1 +
2\alpha ,2\beta ;\xi ,\eta }\right)\displaystyle \frac{d}{ds} \left\{y\left( s \right)\right\}ds \\
& \hskip 6mm + 2\left( {\alpha  + \beta } \right)k_1 y_0 \displaystyle \mathop
{\lim } \limits_{\rho  \to 0} \int\limits_{C_\rho  } {} x^{2\alpha
} y^{2\beta } \left( {r^2 } \right)^{ - \alpha  - \beta  - 1} F_2
\left( {1 + \alpha  + \beta,\alpha ,1 + \beta ;2\alpha ,1 +
2\beta ;\xi ,\eta } \right) \displaystyle \frac{d}{ds} \left\{x\left( s \right)\right\}ds \\
& \hskip 3mm:=  - \left( {\alpha  + \beta } \right)k_1 \mathop {\lim }
\limits_{\rho  \to 0} J_1 \left( {x_0 ,y_0 } \right) - 2\left(
{\alpha  + \beta } \right)k_1 x_0 \mathop {\lim }\limits_{\rho \to
0} J_2 \left( {x_0 ,y_0 } \right) + 2\left( {\alpha  + \beta }
\right)k_1 y_0 \displaystyle \mathop {\lim }\limits_{\rho  \to 0} J_3 \left( {x_0 ,y_0 } \right),
\endaligned \eqno(3.11)
$$
where $J_1$, $J_2$ and $J_3$ are the corresponding integrals in the first equality.
Now, by introducing the polar coordinates:
$$x = x_0  + \rho \cos \varphi \quad \text{and} \quad y =
y_0  + \rho \sin \varphi,$$ we get
$$
\begin{array}{l}
J_1 \left( {x_0 ,y_0 } \right) =  - 2\left( {\alpha  + \beta }
\right)k_1 \displaystyle \mathop {\lim }\limits_{\rho  \to 0}
\displaystyle \int\limits_0^{2\pi } {} \left( {x_0  + \rho \cos \varphi }
\right)^{2\alpha } \left( {y_0  + \rho \sin \varphi } \right)^{2\beta }  \\
\hskip 20mm
\cdot \left( {\rho ^2 } \right)^{ - \alpha  - \beta }
\displaystyle F_2 \left( {\alpha  + \beta +1, \alpha ,
\beta ;2\alpha ,2\beta ;\xi ,\eta ,\zeta } \right)d\varphi . \\
\end{array}\eqno
(3.12)
$$
By using the following known formulas (see \cite[p. 253, Eq.~(26)]{Bu-Ch-40}; also see \cite[p. 113, Eq.~(4)]{Er-Ma-Ob-Tr}):
$$
\begin{array}{l}
\displaystyle F_2 \left( {a,b_1 ,b_2 ;c_1 ,c_2 ;x,y} \right)
= \displaystyle \sum\limits_{j = 0}^\infty  {} \displaystyle
\frac{{\left( a \right)_j \left( {b_1 } \right)_j \left( {b_2 }
\right)_j }}{{\left( {c_1 } \right)_j \left( {c_2 } \right)_j
j!}}x^j y^j\\
\hskip 10mm \cdot {}_2F_1\left( {a + j,b_1  + j;c_1  + j;x} \right){}_2F_1\left( {a
+ j,b_2  + j;c_2  + j;y} \right)
\end{array}\eqno (3.13)
$$
and
$$
\displaystyle {}_2F_1\left( {a,b;c,x} \right) = \left( {1 - x} \right)^{
- b} {}_2F_1\left( {c - a,b;c, \displaystyle \frac{x}{{x - 1}}}
\right),\eqno (3.14)
$$
we obtain
$$
\begin{array}{l}
\displaystyle F_2 \left( {a,b_1 ,b_2 ;c_1 ,c_2 ;x,y} \right)\\
 \hskip 3mm = \left( {1 - x} \right)^{ - b_1 } \left( {1 - y} \right)^{ - b_2
} \displaystyle \sum\limits_{j = 0}^\infty  {}
\displaystyle \frac{{\left( a \right)_j \left( {b_1 }
\right)_j \left( {b_2 } \right)_j }}{{\left( {c_1 } \right)_j \left( {c_2 }
\right)_j j!}}\left( {\displaystyle \frac{x}{{1 - x}}}
\right)^j \left( {\displaystyle \frac{y}{{1 - y}}} \right)^j   \\
 \hskip 5mm\cdot {}_2F_1\left( {c_1  - a,b_1  + j;c_1  + j; \displaystyle
\frac{x}{{x - 1}}} \right){}_2F_1\left( {c_2  - a,b_2  + j;c_2  + j; \displaystyle \frac{y}{{y - 1}}} \right), \\
\end{array}\eqno
(3.15)
$$
where ${}_2F_1\left( {a,b;c;x} \right)$ is Gauss's hypergeometric function
 (see \cite[p. 69,  Eq.~(2)]{Er-Ma-Ob-Tr}). Hence we have
$$
\begin{array}{l}
\displaystyle F_2 \left( {1 + \alpha  + \beta ;\alpha ,\beta ;2\alpha ,2\beta ;\xi ,\eta } \right) \\
\hskip 3mm= \left( {\rho ^2 } \right)^{\alpha  + \beta }
\left( {\rho ^2  + 4x_0^2  + 4x_0 \rho \cos \,\varphi } \right)^{ - \alpha }
\left( {\rho ^2  + 4y_0^2  + 4y_0 \rho \sin \,\varphi } \right)^{ - \beta } P_{11} , \\
\end{array}\eqno
(3.16)
$$
where
$$
\begin{array}{l}
P_{11}  = \displaystyle \sum\limits_{j = 0}^\infty  {}
\displaystyle \frac{{\left( {1 + \alpha  + \beta } \right)_j
\left( \alpha \right)_j \left( \beta  \right)_j }}{{\left(
{2\alpha } \right)_j \left( {2\beta } \right)_j j!j!}}\left( {
\displaystyle \frac{{4x_0^2 + 4x_0 \rho \cos \,\varphi }}{{\rho ^2
+ 4x_0^2  + 4x_0 \rho \cos \,\varphi }}} \right)^j \left(
{\displaystyle \frac{{4y_0^2  + 4y_0 \rho \sin \,
\varphi }}{{\rho ^2  + 4y_0^2  + 4y_0 \rho \sin \,\varphi }}} \right)^j  \\
\,\,\,\,\,\,\,\,\,\,\,\,\,\,\, \cdot {}_2F_1\left( {\alpha  - \beta  -
1,\alpha + j;2\alpha + j;
\displaystyle \frac{{4x_0^2  + 4x_0 \rho \cos \,\varphi }}{{\rho ^2  + 4x_0^2  + 4x_0 \rho \cos \,\varphi }}} \right) \\
\,\,\,\,\,\,\,\,\,\,\,\,\,\,\, \cdot {}_2F_1\left( {\beta  - \alpha  -
1,\beta + j;2\beta + j;\displaystyle \frac{{4y_0^2  + 4y_0
\rho \sin \,\varphi }}{{\rho ^2  + 4y_0^2  + 4y_0 \rho \sin \,\varphi }}} \right). \\
\end{array}
$$
Using the well-known Gauss's summation formula for ${}_2F_1$  (see \cite[p. 112,  Eq.~(46)]{Er-Ma-Ob-Tr})
$$ {}_2F_1\left( {a,b;c;1} \right)  = \displaystyle \frac{{\Gamma \left( c
\right)\Gamma \left( {c - a - b} \right)}}{{\Gamma \left( {c - a}
\right)\Gamma \left( {c - b} \right)}} \quad  \left(\Re \left( {c - a - b} \right) > 0; \,\,c \ne 0, - 1, -
2,...\right), $$
we obtain
$$
\displaystyle \mathop {\lim }\limits_{\rho  \to 0} P_{11}  =
\displaystyle \frac{{\Gamma \left( {2\alpha } \right)\Gamma \left(
{2\beta } \right)}}{{\Gamma \left( \alpha \right)\Gamma \left(
\beta \right)\Gamma \left( {1 + \alpha  + \beta } \right)}}.\eqno
(3.17)
$$
Thus,   by virtue of the identities (3.12), (3.16), and (3.17),
we get
$$
\displaystyle \left( {\alpha  + \beta } \right)k_1 \displaystyle
\mathop {\lim } \limits_{\rho \to 0} J_1 \left( {x_0 ,y_0 }
\right) = 1.\eqno (3.18)
$$
Similarly, by considering  the corresponding identities  and the fact that
   $$\mathop {\lim }\limits_{\rho
\to 0} \rho \ln \rho  = 0, $$ we find that
$$
\displaystyle 2\left( {\alpha  + \beta } \right)k_1 x_0
\displaystyle \mathop {\lim }\limits_{\rho  \to 0} J_2 \left( {x_0
,y_0 } \right) = 2\left( {\alpha  + \beta } \right)k_1 y_0
\displaystyle \mathop {\lim }\limits_{\rho \to 0} J_3 \left( {x_0
,y_0 } \right) = 0.\eqno (3.19)
$$
Hence, by view of (3.18) and (3.19),  the formula  (3.11) in the case of $\left( {x_0
,y_0 } \right) \in \Omega $  becomes
$$
w_1^{\left( 1 \right)} \left( {x_0 ,y_0 } \right) =  - 1.\eqno
(3.20)
$$

\vskip 3mm
\textbf{Case 2.}  When $\left( {x_0 ,y_0 } \right) \in \Gamma$,
 we  cut a circle $C_\rho$  centered at $\left( {x_0 ,y_0 } \right)$  with  a small radius
$\rho$  off  the domain $\Omega$
and denote the remaining part of the curve  by $\Gamma- \Gamma _\rho$.
Let $C_\rho^{'}$   denote  a part of the circle $C_\rho$
lying inside  the domain $\Omega $.
 We  consider the domain
$\Omega _\rho$ which is bounded by a curve $\Gamma - \Gamma _\rho
$, $C_\rho ^{'}$ and the segments $\left[ {0,a} \right]$ and  $\left[
{0,b} \right]$ along the $x-$ and $y-$axes, respectively. Then we have
$$
\displaystyle w_1^{\left( 1 \right)} \left( {x_0 ,y_0 } \right) =
\displaystyle \int_0^l\, {x^{2\alpha } y^{2\beta } }
\displaystyle \frac{\partial}{\partial n} \left\{q_1 \left( {x,y;x_0 ,y_0 }
\right)\right\}ds = \displaystyle \mathop {\lim
}\limits_{\rho  \to 0} \int\limits_{\Gamma - \Gamma _\rho  }
{x^{2\alpha }} y^{2\beta } \displaystyle \frac{\partial}{\partial n} \left\{q_1 \left( {x,y;x_0 ,y_0 }
\right)\right\}ds.\eqno (3.21)
$$
When the point $\left( {x_0 ,y_0 } \right)$ lies outside the domain
$\Omega _\rho$, it is found that, in this domain, $q_1 \left( {x,y;x_0 ,y_0 }
\right)$ is a regular solution of the equation $\left(H_{\alpha ,\beta
}^0\right)$. Therefore, by virtue of (2.6),  we have
$$
\displaystyle \int\limits_{\Gamma  - \Gamma _\rho  } {} x^{2\alpha
} y^{2\beta }\displaystyle \frac{\partial }{{\partial n}} \left\{q_1
\left( {x,y;x_0 ,y_0 } \right)\right\}ds = \displaystyle
\int\limits_{C_\rho ^{'} } {} x^{2\alpha } y^{2\beta }
\displaystyle \frac{\partial }{{\partial n}}\left\{q_1 \left( {x,y;x_0
,y_0 } \right)\right\}ds.\eqno (3.22)
$$
Substituting from (3.22) into (3.21), we get
$$
\displaystyle w_1^{\left( 1 \right)} \left( {x_0 ,y_0 } \right) =
\displaystyle \int_0^l \,{x^{2\alpha } y^{2\beta } }
\displaystyle \frac{\partial}{\partial n} \left\{q_1 \left( {x,y;x_0 ,y_0 }
\right)\right\}ds = \displaystyle \mathop {\lim
}\limits_{\rho  \to 0} \int\limits_{C_\rho ^{'} } {} x^{2\alpha }
y^{2\beta } \displaystyle \frac{\partial}{\partial n} \left\{ q_1 \left( {x,y;x_0 ,y_0
} \right)\right\} ds.\eqno (3.23)
$$
Similarly, by again introducing  the polar coordinates centered at
the point $\left( {x_0 ,y_0 } \right)$, we find that
$$
w_1^{\left( 1 \right)} \left( {x_0 ,y_0 } \right) =  -
\displaystyle \frac{1}{2}.\eqno (3.24)
$$

\vskip 3mm
\textbf{Case 3.} When $\left( {x_0 ,y_0 } \right) \not\in \bar \Omega$,
It is noted that the function $q_1 \left( {x,y;x_0 ,y_0 } \right)$ is a
regular solution of the equation $\left(H_{\alpha ,\beta }^0\right)$. Hence, in view of
 the formula (2.6), we have
$$
\displaystyle w_1^{\left( 1 \right)} \left( {x_0 ,y_0 } \right) =
\displaystyle \int_0^l \, {x^{2\alpha } y^{2\beta } }
\displaystyle \frac{\partial }{{\partial n}} \left\{q_1 \left( {x,y;x_0
,y_0 } \right)\right\} ds = 0.\eqno (3.25)
$$
The proof of Lemma 1 is thus complete.
\end{proof}

\vskip 5mm

\textbf{Lemma 2.} \emph{The following formula holds true}:
$$
w_1^{\left( 1 \right)} \left( {x_0 ,0 } \right) =
\left\{
  \begin{array}{rl}
  - 1 & \qquad  \left(x_0  \in \left( {0,a} \right)\right) \\
    - \frac{1}{2} &  \qquad \left( x_0  = 0\,\,\,\text{or}\,\,\,x_0  = a\right)\\
    0 &  \qquad \left( a < x_0 \right).
  \end{array}
\right.\eqno (3.26)
$$

\begin{proof} \textbf{Case 1.} When $x_0  \in \left(0,a\right)$, we introduce
a straight line $y  = h$ for a sufficiently small positive real number $h$  and
consider a domain $\Omega _h $ which is the part of the domain $\Omega $
lying above the straight line $y  = h$. Applying the formula
(2.6), we obtain
$$
\displaystyle w_1^{\left( 1 \right)} \left( {x_0 ,0} \right) =
\displaystyle \mathop {\lim }\limits_{h \to 0} \displaystyle
\int_0^{x_1 } {} \left. {x^{2\alpha } y^{2\beta }
\displaystyle \frac{{\partial q_1 \left( {x,y;x_0 ,0}
\right)}}{{\partial y}}} \right|_{y = h} dx,\eqno (3.27)
$$
where $x_1 \left( \varepsilon  \right)$ is an abscissa of a point at which
the straight line $y  = h$ intersects the curve $\Gamma $.
It follows from (3.8) and (3.27) that
$$
\displaystyle w_1^{\left( 1 \right)} \left( {x_0 ,0} \right) =  -
2\left( {\alpha  + \beta } \right)k_1 \displaystyle \mathop {\lim
}\limits_{h \to 0} h^{1 + 2\beta } \int_0^{x_1 } {}
x^{2\alpha } \displaystyle \frac{{F\left( {\alpha  + \beta  +
1,\alpha ;2\alpha ,\displaystyle \frac{{ - 4xx_0 }}{{\left( {x -
x_0 } \right)^2 + h^2 }}} \right)}}{{\left[ {\left( {x - x_0 }
\right)^2  + h^2 } \right]^{\alpha  + \beta  + 1} }}dx.\eqno
(3.28)
$$
Using the hypergeometric transformation formula (3.14) inside the integrand  of (3.28), we have
$$
\displaystyle w_1^{\left( 1 \right)} \left( {x_0 ,0} \right) =  -
2\left( {\alpha  + \beta } \right)k_1 \displaystyle \mathop {\lim
}\limits_{h \to 0} h^{1 + 2\beta } \int_0^{x_1 } {}
x^{2\alpha } \displaystyle \frac{{F\left( {\alpha  - \beta  -
1,\alpha ;2\alpha ,\displaystyle \frac{{4xx_0 }}{{\left( {x + x_0
} \right)^2  + h^2 }}} \right)}}{{\left[ {\left( {x - x_0 }
\right)^2  + h^2 } \right]^{\beta  + 1} \left[ {\left( {x + x_0 }
\right)^2  + h^2 } \right]^\alpha  }}dx,\eqno (3.29)
$$
which, upon setting $x = x_0  + ht$ inside the integrand,
yields
$$
\displaystyle w_1^{\left( 1 \right)} \left( {x_0 ,0} \right) =  -
2\left( {\alpha  + \beta } \right)k_1 \displaystyle \mathop {\lim
}\limits_{h \to 0} \displaystyle \int_{l_1 }^{l_2 } {}
\left( {x_0  + ht} \right)^{2\alpha } \displaystyle \frac{{F\left(
{\alpha  - \beta  - 1,\alpha ;2\alpha ,\displaystyle \frac{{4x_0
\left( {x_0  + ht} \right)}}{{\left( {2x_0 + ht} \right)^2  + h^2
}}} \right)}}{{\left( {1 + t^2 } \right)^{\beta  + 1} \left[
{\left( {2x_0  + ht} \right)^2  + h^2 } \right]^\alpha  }}dt,\eqno
(3.30)
$$
where
$$l_1  =  -  \frac{{x_0 }}{h} \quad  \text{and} \quad   l_2  =
 \frac{{x_1 - x_0 }}{h}.
$$
Considering
$$
\displaystyle \mathop {\lim }\limits_{h \to 0} F\left( {\alpha  -
\beta  - 1,\alpha ;2\alpha ,\displaystyle \frac{{4x_0 \left( {x_0
+ ht} \right)}}{{\left( {2x_0  + ht} \right)^2  + h^2 }}} \right)
= F\left( {\alpha  - \beta  - 1,\alpha ;2\alpha ,1} \right) =
\displaystyle \frac{{\Gamma \left( {2\alpha } \right)\Gamma \left(
{1 + \beta } \right)}}{{\left( {\alpha  + \beta } \right)\Gamma
\left( {\alpha + \beta } \right)\Gamma \left( \alpha  \right)}}
$$
and
$$
\displaystyle \int_{ - \infty }^{ + \infty } {}
\displaystyle \frac{{dt}}{{\left( {1 + t^2 } \right)^{\beta  + 1}
}} = \displaystyle \frac{{\pi \Gamma \left( {2\beta }
\right)}}{{2^{2\beta  - 1} \beta \Gamma ^2 \left( \beta \right)}},
$$
we find from (3.30) that
$$
\displaystyle w_1^{\left( 1 \right)} \left( {x_0 ,0} \right) =  -
1.\eqno (3.31)
$$

\vskip 3mm
The other cases when $x_0=0$, $x_0=a$ and $x_0>a$ will be proved by using
the similar argument as in Case 1.

This evidently completes our proof of Lemma 2.
\end{proof}
\vskip 5mm

\textbf{Lemma 3.} \emph{The following formula holds true}:
$$
w_1^{\left( 1 \right)} \left( {0 ,y_0 } \right) =
\left\{
  \begin{array}{rl}
  - 1 &  \qquad  \left( y_0  \in \left( {0,b} \right)\right) \\
    - \frac{1}{2} &  \qquad  \left( y_0  = 0\,\,\text{or}\,\,y_0  = b\right)\\
    0 &  \qquad  \left( b < y_0\right).
  \end{array}
\right.\eqno (3.32)
$$

\begin{proof}
The proof of Lemma 3 would run parallel to that of Lemma 2.
\end{proof}

\vskip 5mm
\textbf{Theorem 1.} \emph{For any points} $\left( {x,y}
\right)$ \emph{and}  $\left( {x_0 ,y_0 } \right) \in \mathbb{R}_ + ^2 $ and $x \ne
x_0$ and $y \ne y_0$,  \emph{the following inequality holds true}:
$$
\left| {q_1 \left( {x,y;x_0 ,y_0 } \right)} \right| \le k_1
\displaystyle \frac{{\Gamma \left( {2\alpha } \right)\Gamma \left(
{2\beta } \right)}}{{\Gamma ^2 \left( {\alpha  + \beta }
\right)}}\left( {r_1^2 } \right)^{ - \alpha } \left( {r_2^2 }
\right)^{ - \beta } F\left[ {\alpha ,\beta ;\alpha  + \beta
;\left( {1 - \displaystyle \frac{{r^2 }}{{r_1^2 }}} \right)\left(
{1 - \displaystyle \frac{{r^2 }}{{r_2^2 }}} \right)}
\right],\eqno(3.33)
$$
\emph{where} $\alpha$ \emph{and} $\beta$ \emph{are real parameters with} $\left(0<\alpha,\,\beta<\frac{1}{2}\right)$ \emph{as in the equation}
$\left( {H_{\alpha ,\beta }^\lambda  } \right)$, \emph{and} $r$, $r_1$ \emph{and} $r_2$ \emph{are as in} (1.10).

\begin{proof}  It follows from (3.15) that
$$
\begin{array}{l}
q_1 \left( {x,y;x_0 ,y_0 } \right) = k_1 \left( {r_1^2 } \right)^{
- \alpha } \left( {r_2^2 } \right)^{ - \beta } \displaystyle
\sum\limits_{j = 0}^\infty  {} \displaystyle \frac{{\left( {\alpha
+ \beta } \right)_j \left( \alpha  \right)_j \left( \beta
\right)_j }}{{\left( {2\alpha } \right)_j \left( {2\beta }
\right)_j j!}}\left( {1 - \displaystyle \frac{{r^2 }}{{r_1^2 }}}
\right)^j \left( {1 - \displaystyle \frac{{r^2 }}{{r_2^2 }}} \right)^j  \\
 \hskip 13mm \cdot {}_2F_1\left( {\alpha  - \beta ,\alpha  + j;2\alpha  + j;1 -
\displaystyle \frac{{r^2 }}{{r_1^2 }}}
\right){}_2F_1\left( {\beta  - \alpha ,\beta  + j;2\beta  + j;1 - \displaystyle \frac{{r^2 }}{{r_2^2 }}} \right). \\
\end{array}\eqno (3.34)
$$
Now, in view of the following inequalities:
$$
{}_2F_1\left( {\alpha  - \beta ,\alpha  + j;2\alpha  + j;1 -
\displaystyle \frac{{r^2 }}{{r_1^2 }}} \right) \le \displaystyle
\frac{{\Gamma \left( {2\alpha } \right)\Gamma \left( \beta
\right)\left( {2\alpha } \right)_j }}{{\Gamma \left( {\alpha  +
\beta } \right)\Gamma \left( \alpha \right)\left( {\alpha  + \beta
} \right)_j }}
$$
and
$$
{}_2F_1\left( {\beta  - \alpha ,\beta  + j;2\beta  + j;1 - \displaystyle
\frac{{r^2 }}{{r_2^2 }}} \right) \le \displaystyle \frac{{\Gamma
\left( {2\beta } \right)\Gamma \left( \alpha  \right)\left(
{2\beta } \right)_j }}{{\Gamma \left( {\alpha  + \beta }
\right)\Gamma \left( \beta \right)\left( {\alpha  + \beta }
\right)_j }},
$$
we find from (3.34) that the inequality (3.33) holds.
Hence Theorem 1 is proved.
\end{proof}

\vskip 5 mm

 By virtue of the following  known formula
\cite[p. 117,  Eq.~(12)]{Er-Ma-Ob-Tr}
$$
\aligned
{}_2F_1\left( {a,b;a + b;z} \right)  = & - \displaystyle \frac{{\Gamma
\left( {a + b} \right)}}{{\Gamma \left( a \right)\Gamma \left( b
\right)}}F\left({a,b;1;1 - z} \right)\ln \left( {1 - z} \right) \\
&\hskip -10mm  + \displaystyle \frac{{\Gamma \left( {a + b} \right)}}{{\Gamma ^2
\left( a \right)\Gamma ^2 \left( b \right)}} \displaystyle
\sum\limits_{j = 0}^\infty  {} \displaystyle \frac{{\Gamma \left(
{a + j} \right)\Gamma \left( {b + j} \right)}}{{\left( {j!}
\right)^2 }}\left[ {2\psi \left( {1 + j} \right) - \psi \left( {a
+ j} \right) - \psi \left( {b + j} \right)} \right]\left( {1 - z} \right)^j \\
& \left(- \pi  < \arg \,\left( {1 - z} \right) < \pi;\,\,a,b \ne 0, - 1,
- 2, \cdots\right),
\endaligned
$$
where
  $$ \psi (z):= \frac{d}{dz} \left\{\ln \Gamma (z)\right\}=\frac{\Gamma' (z)}{\Gamma (z)}
\quad \text{or} \quad \ln \Gamma (z)= \int_1^z\, \psi (t)\,dt, $$
we observe from (3.33) that $ q_1 \left( {x,y;x_0 ,y_0 } \right)$ has a
logarithmic singularity at $r = 0$.

\vskip 5 mm
\textbf{Theorem 2.} \emph{If the curve} $ \Gamma$ \emph{satisfies
conditions} (3.1), \emph{then the following inequality holds true}:
$$
\displaystyle \int\limits_\Gamma  {} x^{2\alpha } y^{2\beta }
\left| {\displaystyle \frac{{\partial q_1 \left( {x,y;x_0 ,y_0 }
\right)}}{{\partial n}}} \right|ds \le C_1,
$$
\emph{where} $ C_1$ \emph{is a constant}.

\begin{proof}
 Theorem 2 follows by suitably applying Lemmas 1 to 3.
\end{proof}

\vskip 5 mm
\textbf{Theorem 3.} \emph{If} $\mu _1 \left( t \right) \in \left[
{0,\,\,l} \right]$,  \emph{then the following limiting formulas hold true for a double-layer potential} (3.2):
$$
w_i^{\left( 1 \right)} \left( t \right) =  - \displaystyle
\frac{1}{2}\mu _1 \left( t \right) + \displaystyle \int_0^l
{} \mu _1 \left( s \right)K_1 \left( {s,t} \right)ds
$$
\emph{and}
$$
w_e^{\left( 1 \right)} \left( t \right) =\displaystyle
\frac{1}{2}\mu _1 \left( t \right) + \displaystyle \int_0^l
{} \mu _1 \left( s \right)K_1 \left( {s,t} \right)ds,
$$
\emph{where}
$$K_1 \left( {s,t} \right) = x^{2\alpha } \left( s
\right)y^{2\beta } \left( s \right)\, \frac{\partial}{\partial n} \left\{q_1 \left[ {x\left( s
\right),y\left( s \right);x_0 \left( t \right),y_0 \left( t
\right)} \right]\right\}$$
$$ \Big(\left( {x\left( s \right),y\left( s \right)} \right)\in \Gamma; \,\, \left( {x_0 \left( t
\right),y_0 \left( t \right)} \right) \in \Gamma \Big),
$$
 $w_i^{\left( 1 \right)} \left( t \right)$ \emph{and}  $w_e^{\left( 1
\right)} \left( t \right)$  \emph{are limiting values of the double-layer
potential}  (3.2) \emph{at}
$$\left( {x_0 \left( t \right),y_0 \left( t \right)} \right) \to \Gamma $$
\emph{from the inside and the outside}, \emph{respectively}.

\begin{proof}
 We find from Lemma 1, in conjunction with Theorems 1 and 2, that each of the limiting formulas
 asserted by Theorem 3 holds true.
\end{proof}

\vskip 2mm

\begin{center}
\textbf{Acknowledgements}\vskip 1mm
\end{center}

This research was supported by the {\it Basic Science Research Program} through the {\it National Research Foundation of Korea} funded by the Ministry of Education, Science and Technology of the Republic of Korea (2011-0005224). The present investigation was supported, in part, by the {\it Natural Sciences and Engineering Research Council of Canada} under Grant OGP0007353.

\vskip 10mm

\bibliographystyle{amsplain}

\begin{thebibliography}{99}

\bibitem{Alti-1}
A. Altin,  Solutions of type $r^m $ for a class of singular
equations, \emph{Internat. J.  Math. Sci.}
\textbf{5} (1982), 613--619.




\bibitem{Alti-2}
A. Altin, Some expansion formulas for a class of singular
partial differential equations, \emph{Proc. Amer.
Math. Soc.} \textbf{85} (1982), 42--46.


\bibitem{Al-Eu}
A. Altin and Y. Eutiquio,  Some properties of solutions of a
class of singular partial differential equations, \emph{Bull.
 Inst. Math.  Acad. Sinica} \textbf{11} (1983), 81--87.

\bibitem{Ap-Ka}
P. Appell and J. Kamp\'e de F\'eriet, \emph{Fonctions
Hyperg\'eom\'etriques et Hypersph\'eriques; Polyn$\hat{o}$mes d'Hermite},
Gauthier-Villars, Paris, 1926.

\bibitem{Bu-Ch-40}
J. L. Burchnall and T. W. Chaundy, Expansions of Appell's double
hypergeometric functions, \emph{Quart. J. Math. Oxford Ser.} \textbf{11}
(1940), 249--270.


\bibitem{Erde-56}
A. Erd\'{e}lyi, Singularities of generalized axially symmetric
potentials, \emph{Comm. Pure Appl. Math.} \textbf{2} (1956), 403--414.



\bibitem{Erde-65}
A. Erd\'{e}lyi,  An application of fractional integrals, \emph{J.
Analyse Math.} \textbf{14} (1965), 113--126.


\bibitem{Er-Ma-Ob-Tr}
A.  Erd\'elyi, W.  Magnus,  F. Oberhettinger and F. G. Tricomi,
{\it  Higher Transcendental Functions}, Vol. {\bf I},
McGraw-Hill Book Company, New York, Toronto and London, 1953;
Russian edition, Izdat. Nauka, Moscow,  1973.




\bibitem{Frya}
A. J. Fryant,  Growth and complete sequences of generalized
bi-axially symmetric potentials, \emph{J. Differential
Equations} \textbf{31} (1979), 155--164.




\bibitem{Gilb-60}
R. P. Gilbert,  On the singularities of generalized axially
symmetric potentials, \emph{Arch. Rational Mech. Anal.} \textbf{6} (1960),
171--176.

\bibitem{Gilb-62}
R. P. Gilbert,  Some properties of generalized axially symmetric
potentials, \emph{Amer. J. Math.} \textbf{84} (1962), 475--484.

\bibitem{Gilb-64}
R. P. Gilbert,  Bergman's integral operator method in
generalized axially symmetric potential theory, \emph{J.
Math. Phys.} \textbf{5} (1964), 983--987.


\bibitem{Gilb-65}
R. P. Gilbert,  On the location of singularities of a class of
elliptic partial differential equations in four variables,
\emph{Canad. J. Math.} \textbf{17} (1965), 676--686.

\bibitem{Gilb-67}
R. P. Gilbert,  On the analytic properties of solutions to a
generalized axially symmetric Schr$\ddot{\rm o}$dinger equations,
\emph{J. Differential Equations} \textbf{3} (1967), 59--77.


\bibitem{Gilb-68}
R. P. Gilbert,  An investigation of the analytic properties of
solutions to the generalized axially symmetric, reduced wave
equation in $n + 1$ variables, with an application to the theory
of potential scattering, \emph{SIAM J. Appl. Math.} \textbf{16} (1968), 30--50.


\bibitem{Gilb-69}
R. P. Gilbert,  \emph{Function Theoretic Methods in Partial
Differential Equations}, Mathematics in Science and Engineering, Vol. \textbf{54},
A Series of Monographs and Textbooks, Academic Press, New York and London, 1969.




\bibitem{Gi-Ho}
R. P. Gilbert and H. Howard,  On solutions of the generalized
axially symmetric wave equation represented by Bergman operators,
\emph{Proc. London Math. Soc.} (\emph{Ser}. 3) \textbf{15} (1965), 346--360.

\bibitem{gunt-67}
 N. M. G$\ddot{\rm u}$nter,   \emph{Potential Theory and Its Applications to Basic Problems of Mathematical Physics}
 (Translated from the Russian edition by J. R. Schulenberger), Frederick Ungar Publishing Company,
 New York, 1967.


\bibitem{Hasa-07}
A. Hasanov,  Fundamental solutions of generalized bi-axially
symmetric Helmholtz equation,  \emph{Complex Variables and
Elliptic Equations} \textbf{52} (2007), 673--683.



\bibitem{Henr-53}
P. Henrici,  Zur Funktionentheorie der Wellengleichung,
\emph{Comment. Math. Helv.} \textbf{27} (1953), 235--293.


\bibitem{Henr-57}
P. Henrici,  On the domain of regularity of generalized
axially symmetric potentials. \emph{Proc. Amer. Math. Soc.} \textbf{8} (1957),
29--31.


\bibitem{Henr-60}
P. Henrici,  Complete systems of solutions for a class
of singular elliptic partial differential equations, in  \emph{Boundary
Value Problems in Differential Equations}, pp. 19--34, University of Wisconsin
Press, Madison,  1960.


\bibitem{Hube}
A. Huber,  On the uniqueness of generalized axisymmetric
potentials,  \emph{Ann. Math.} \textbf{60} (1954), 351--358.



\bibitem{Kuma}
 D. Kumar, Approximation of growth numbers generalized
bi-axially symmetric potentials, \emph{Fasc. Math.}, \textbf{35} (2005),
51--60.


\bibitem{Lo}
 C. Y.  Lo,  Boundary value problems of generalized axially
symmetric Helmholtz equations, \emph{Portugal. Math.}
\textbf{36} (1977), 279--289.


\bibitem{Mari}
O. I. Mari$\breve{\rm c}$ev, An integral representation of the solutions of a
generalized biaxially symmetric Helmholtz equation and formulas for its inversion  (Russian),
\emph{Differencial'nye Uravnenija} \textbf{14} (1978), 1824--1831.



\bibitem{McCo-79}
P. A. McCoy,  Polynomial approximation and growth of
generalized axisymmetric potentials, \emph{Canad. J.
Math.} \textbf{31} (1979), 49--59.



\bibitem{McCo-80}
P. A. McCoy,  Best $L^p $-approximation of generalized
bi-axisymmetric potentials,  \emph{Proc. Amer.
Math. Soc.}  \textbf{79} (1980),  435--440.

\bibitem{miran-70}
C. Miranda,   \emph{Partial Differential Equations of Elliptic Type},
Second Revised Edition (Translated from the Italian by Z. C. Motteler),
Ergebnisse der Mathematik und ihrer Grenzgebiete, Band \textbf{2},
Springer-Verlag, Berlin, Heidelberg and New York,    1970.


\bibitem{Pi-Bo}
P. P. Niu  and X. B. Lo,  Some notes on solvability of LPDO,
\emph{J. Math. Res. Exposition} \textbf{3} (1983),
127--129.



\bibitem{Rang}
K. B. Ranger,  On the construction of some integral operators
for generalized axially symmetric harmonics and stream functions,
\emph{J. Math. Mech.} \textbf{14} (1965), 383--401.


\bibitem{Ra-Ha}
J. M. Rassias  and A. Hasanov,   Fundamental solutions of two
degenerated elliptic equations and solutions of boundary value
problems in infinite area, \emph{Internat. J. Appl.
Math. Statist.} \textbf{8} (M07) (2007), 87--95.

\bibitem{Sr-Ka}
H. M. Srivastava and P. W. Karlsson, \emph{Multiple Gaussian
Hypergeometric Series}, Halsted Press (Ellis Horwood Limited,
Chichester), John Wiley and Sons, New York, Chichester, Brisbane  and Toronto,
1985.


\bibitem{Wein}
R. J. Weinacht,  Some properties of generalized axially
symmetric Helmholtz potentials, \emph{SIAM J. Math. Anal}. \textbf{5} (1974),
147--152.



\bibitem{Wein-48}
A. Weinstein,  Discontinuous integrals and generalized
potential theory, \emph{Trans. Amer. Math. Soc.} \textbf{63} (1948), 342--354.



\bibitem{Wein-53}
A. Weinstein,  Generalized axially symmetric potential
theory, \emph{Bull. Amer. Math. Soc.} \textbf{59} (1953), 20--38.













\end{thebibliography}

\end{document}